\documentclass[10pt]{article}
\usepackage{cite}
\usepackage{mathrsfs}
\usepackage{amsfonts}
\usepackage{amsmath}
\usepackage{amsfonts,amssymb,color}
\usepackage{dsfont}
\usepackage{curves}
\usepackage{mathrsfs}
\usepackage{pifont}
\usepackage{amssymb}
\allowdisplaybreaks

\numberwithin{equation}{section}

\date{}

\def\BigRoman{\uppercase\expandafter{\romannumeral\number\count 255 }}
\def\Romannumeral{\afterassignment\BigRoman\count255=}

\begin{document}
\title{Adjacency spectral radius and $H$-factors in 1-binding graphs
}
\author{\small  Sizhong Zhou$^{1}$\footnote{Corresponding
author. E-mail address: zsz\_cumt@163.com (S. Zhou)}, Tao Zhang$^{2}$, Zhiren Sun$^{3}$\\
\small $1$. School of Science, Jiangsu University of Science and Technology,\\
\small Zhenjiang, Jiangsu 212100, China\\
\small $2$. School of Economics and management, Jiangsu University of Science and Technology,\\
\small Zhenjiang, Jiangsu 212100, China\\
\small $3$. School of Mathematical Sciences, Nanjing Normal University,\\
\small Nanjing, Jiangsu 210023, China\\
}

\maketitle
\begin{abstract}
\noindent Let $G$ be a graph, and let $H:V(G)\longrightarrow\{\{1\},\{0,2\}\}$ be a set-valued function. Hence, $H(v)$ equals $\{1\}$ or $\{0,2\}$ for any $v\in V(G)$. We let
$$
H^{-1}(1)=\{v: v\in V(G) \ \mbox{and} \ H(v)=1\}.
$$
An $H$-factor of $G$ is a spanning subgraph $F$ of $G$ such that $d_F(v)\in H(v)$ for each $v\in V(G)$. Lu and Kano showed a characterization for the existence of an $H$-factor in a graph
[Characterization of 1-tough graphs using factors, Discrete Math. 343 (2020) 111901]. Let $A(G)$ and $\rho(G)$ denote the adjacency matrix and the adjacency spectral radius of $G$, respectively.
By using Lu and Kano's result, we pose a sufficient condition with respect to the adjacency spectral radius to guarantee the existence of an $H$-factor in a 1-binding graph. In this paper, we
prove that if a connected 1-binding graph $G$ of order $n\geq11$ satisfies $\rho(G)\geq\rho(K_1\vee(K_{n-4}\cup K_2\cup K_1))$, then $G$ has an $H$-factor for each $H:V(G)\longrightarrow\{\{1\},\{0,2\}\}$
with $H^{-1}(1)$ even, unless $G=K_1\vee(K_{n-4}\cup K_2\cup K_1)$.
\\
\begin{flushleft}

(2020) Mathematics Subject Classification: 05C50, 05C70
\end{flushleft}
\end{abstract}

\section{Introduction}

Throughout this paper, we only consider finite and undirected graphs without loops or multiple edges. Let $G$ denote a graph with vertex set $V(G)$ and edge set $E(G)$. For $v\in V(G)$, we let
$d_G(v)$ and $N_G(v)$ denote the degree and the neighborhood of $v$ in $G$, respectively. Given $S\subseteq V(G)$, let $N_G(S)$ and $G[S]$ denote the set of neighbors of $S$ in $G$ and the
subgraph of $G$ induced by $S$, respectively. Write $G-S=G[V(G)\setminus S]$. The number of isolated vertices and the number of components in $G$ are denoted by $i(G)$ and $\omega(G)$, respectively.
The complete graph of order $n$ is denoted by $K_n$.

For vertex-disjoint graphs $G_1$ and $G_2$, the graph with vertex set $V(G_1)\cup V(G_2)$ and edge set $E(G_1)\cup E(G_2)$ is called the union of $G_1$ and $G_2$, which is denoted by $G_1\cup G_2$.
For any integer $k\geq1$, let $kG$ denote the disjoint union of $k$ copies of $G$. The join of $G_1$ and $G_2$ is denoted by $G_1\vee G_2$, which is obtained from $G_1\cup G_2$ by joining each
vertex of $G_1$ with each vertex of $G_2$ by an edge.

The binding number was first introduced by Woodall \cite{W} and Anderson \cite{A}. The binding number of $G$, denoted by $bind(G)$, is defined by
$$
bind(G)=\min\left\{\frac{|N_G(S)|}{|S|}:\emptyset\neq S\subseteq V(G), N_G(S)\neq V(G)\right\}.
$$
A graph $G$ is said to be $r$-binding if $bind(G)\geq r$, where $r$ is a positive real number. We refer the reader to \cite{KW,KT,ZS,C,RW,Wr} for some properties on binding numbers.

Let $G$ be a graph with $V(G)=\{v_1,v_2,\ldots,v_n\}$. The adjacency matrix of $G$ is the real symmetric matrix $A(G)=(a_{ij})_{n\times n}$, where $a_{ij}=1$ if $v_iv_j\in E(G)$ and $a_{ij}=0$ if
$v_iv_j\notin E(G)$. The largest eigenvalue of $A(G)$, denoted by $\rho(G)$, is called the adjacency spectral radius of $G$. For some results on spectral radius in graphs, we refer the reader to
\cite{O,Wc,Ws,LM,FLL,Zs,Zt,ZZ,ZZL,ZW}.

Let $H:V(G)\longrightarrow\{\{1\},\{0,2\}\}$ be a set-valued function. Hence, $H(v)$ equals $\{1\}$ or $\{0,2\}$ for any $v\in V(G)$. We let
$$
H^{-1}(1)=\{v: v\in V(G) \ \mbox{and} \ H(v)=1\}.
$$
An $H$-factor of $G$ is a spanning subgraph $F$ of $G$ such that $d_F(v)\in H(v)$ for each $v\in V(G)$. In terms of the Handshaking Lemma, we easily see that $H^{-1}(1)$ is even if $G$ contains
an $H$-factor. Obviously, $G$ contains no $H$-factor if $H^{-1}(1)$ is odd.

Lu and Kano \cite{LK} obtained a characterization for the existence of an $H$-factor in a graph.

\medskip

\noindent{\textbf{Theorem 1.1}} (Lu and Kano \cite{LK}). A connected graph $G$ contains an $H$-factor for each $H:V(G)\longrightarrow\{\{1\},\{0,2\}\}$ with $H^{-1}(1)$ even if and only if
$$
\omega(G-S)\leq|S|+1
$$
for any $S\subseteq V(G)$.

\medskip

Motivated by \cite{LK,FL} directly, we establish a sufficient condition with respect to the adjacency spectral radius to guarantee the existence of an $H$-factor in a 1-binding graph.

\medskip

\noindent{\textbf{Theorem 1.2.}} Let $G$ be a connected 1-binding graph of order $n\geq11$. If
$$
\rho(G)\geq\rho(K_1\vee(K_{n-4}\cup K_2\cup K_1)),
$$
then $G$ has an $H$-factor for each $H:V(G)\longrightarrow\{\{1\},\{0,2\}\}$ with $H^{-1}(1)$ even, unless $G=K_1\vee(K_{n-4}\cup K_2\cup K_1)$.

\medskip

\section{Some preliminaries}

In this section, we provide several necessary preliminary lemmas, which will be used to verify our main result.

\medskip

\noindent{\textbf{Lemma 2.1}} (Li and Feng \cite{LF}). Let $G$ denote a graph, and let $H$ be a subgraph of $G$. Then
$$
\rho(G)\geq\rho(H),
$$
with equality occurring if and only if $G=H$.

\medskip

\noindent{\textbf{Lemma 2.2.}} Let $\sum\limits_{i=1}^{t}n_i=n-s$ with $s\geq1$. If $n_t\geq n_{t-1}\geq\cdots\geq n_1\geq1$, $n_{t-1}\geq p\geq1$ and $n_t<n-s-t-p+2$, then
$$
\rho(K_s\vee(K_{n_t}\cup K_{n_{t-1}}\cup\cdots\cup K_{n_1}))<\rho(K_s\vee(K_{n-s-t-p+2}\cup K_p\cup(t-2)K_1)).
$$

\medskip

\noindent{\it Proof.} Write $G=K_s\vee(K_{n_t}\cup K_{n_{t-1}}\cup\cdots\cup K_{n_1})$, and let $\mathbf{x}$ be the Perron vector of $A(G)$. Using symmetry, we let $\mathbf{x}(v)=x_i$ for any $v\in V(K_{n_i})$,
where $1\leq i\leq t$, and $\mathbf{x}(u)=y$ for any $u\in V(K_s)$. Notice that $K_{n_t+s}$ is a proper subgraph of $G$. Together with Lemma 2.1, we deduce $\rho(G)>\rho(K_{n_t+s})=n_t+s-1>n_t-1$. Notice that
$n_t\geq n_{t-1}\geq\cdots\geq n_1\geq1$. According to $A(G)\mathbf{x}=\rho(G)\mathbf{x}$, we possess
\begin{align*}
(\rho(G)-(n_i-1))(x_t-x_i)=&(\rho(G)-(n_i-1))x_t-(\rho(G)-(n_i-1))x_i\\
=&(\rho(G)-(n_i-1))x_t-sy\\
=&(\rho(G)-(n_t-1)+(n_t-1)-(n_i-1))x_t-sy\\
=&sy+(n_t-n_i)x_t-sy\\
=&(n_t-n_i)x_t\\
\geq&0
\end{align*}
for $1\leq i\leq t$. This leads to $x_t\geq x_i$ for $1\leq i\leq t$. Write $G'=K_s\vee(K_{n-s-t-p+2}\cup K_p\cup(t-2)K_1)$. Then it follows from $p\leq n_{t-1}\leq n_t<n-s-t-p+2$, $n_j\geq1$ for
$1\leq j\leq t-2$ and $x_t\geq x_i$ for $1\leq i\leq t$ that
\begin{align*}
\rho(G')-\rho(G)\geq&\mathbf{x}^{T}(A(G')-A(G))\mathbf{x}\\
=&2(n_{t-1}-p)x_{t-1}(n_tx_t+\sum\limits_{j=1}^{t-2}(n_j-1)x_j-px_{t-1})\\
&+2\sum\limits_{j=1}^{t-2}(n_j-1)x_j(n_tx_t-x_j)\\
&+2\sum\limits_{i=2}^{t-2}\sum\limits_{j=i-1}^{t-3}(n_i-1)(n_j-1)x_ix_j\\
>&0,
\end{align*}
that is, $\rho(K_s\vee(K_{n_t}\cup K_{n_{t-1}}\cup\cdots\cup K_{n_1}))<\rho(K_s\vee(K_{n-s-t-p+2}\cup K_p\cup(t-2)K_1))$. This completes the proof of Lemma 2.2. \hfill $\Box$

\medskip

Let $M$ denote a real matrix of order $n$ and $\mathcal{N}=\{1,2,\ldots,n\}$. For a given partition $\pi:\mathcal{N}=\mathcal{N}_1\cup\mathcal{N}_2\cup\cdots\cup\mathcal{N}_r$, the matrix $M$ can
be correspondingly written as
\begin{align*}
M=\left(
  \begin{array}{cccc}
    M_{11} & M_{12} & \cdots & M_{1r}\\
    M_{21} & M_{22} & \cdots & M_{2r}\\
    \vdots & \vdots & \ddots & \vdots\\
    M_{r1} & M_{r2} & \cdots & M_{rr}\\
  \end{array}
\right),
\end{align*}
where the block $M_{ij}$ denotes the $n_i\times n_j$ matrix for any $1\leq i,j\leq r$. Let $m_{ij}$ denote the average row sum of $M_{ij}$, that is to say, $m_{ij}$ is the sum of all entries in $M_{ij}$
divided by the number of rows. Then the $r\times r$ matrix $M_{\pi}=(m_{ij})$ is called the quotient matrix of $M$. The partition $\pi$ is called equitable if every block $M_{ij}$ of $M$ admits constant
row sum $m_{ij}$ for $1\leq i,j\leq r$.

\medskip

\noindent{\textbf{Lemma 2.3}} (Brouwer and Haemers \cite{BH}, You, Yang, So and Xi \cite{YYSX}). Let $M$ be a real matrix of order $n$ with an equitable partition $\pi$, and let $M_{\pi}$ be the corresponding
quotient matrix. Then the eigenvalues of $M_{\pi}$ are eigenvalues of $M$. Furthermore, if $M$ is nonnegative and irreducible, then the largest eigenvalues of $M$ and $M_{\pi}$ are equal.

\section{Proof of Theorem 1.2}

\noindent{\it Proof of Theorem 1.2.} Suppose that a connected 1-binding graph $G$ has no $H$-factor. Using Theorem 1.1, we conclude $\omega(G-S)\geq|S|+2$ for some nonempty subset $S$ of $V(G)$.
Let $|S|=s$ and $\omega(G-S)=q$. Then $q\geq s+2$. Let $C_1,C_2,\ldots,C_{q-1},C_q$ denote the $q$ components in $G-S$ and $|C_i|=n_i$ for $1\leq i\leq q$. Without loss of generality, we let
$n_q\geq n_{q-1}\geq\cdots\geq n_2\geq n_1\geq1$. We are to verify $n_{s+1}\geq2$. Otherwise, if $n_i=1$ for $1\leq i\leq s+1$, then we let $X=V(C_1)\cup V(C_2)\cup\cdots\cup V(C_{s+1})$.
Obviously, $N_G(X)\subseteq S$ and $|X|=s+1$. Thus, we deduce
$$
\frac{|N_G(X)|}{|X|}\leq\frac{|S|}{|X|}=\frac{s}{s+1}<1,
$$
which is impossible because $G$ is 1-binding. This leads to $n_i\geq2$ for $i\geq s+1$. We easily see that $G$ is a spanning subgraph of $G_1=K_s\vee(K_{n_1}\cup K_{n_2}\cup\cdots\cup K_{n_{s+1}}\cup K_{n_{s+2}})$,
where $n_{s+2}\geq n_{s+1}\geq\cdots\geq n_2\geq n_1\geq1$, $n_{s+1}\geq2$ and $\sum\limits_{i=1}^{s+2}n_i=n-s$. In view of Lemma 2.1, we get
\begin{align}\label{eq:3.1}
\rho(G)\leq\rho(G_1),
\end{align}
with equality occurring if and only if $G=G_1$. Let $G_2=K_s\vee(K_{n-2s-2}\cup K_2\cup sK_1)$, where $n\geq2s+4$. According to Lemma 2.2 (set $p=2$ and $t=s+2$ in Lemma 2.2), we conclude
\begin{align}\label{eq:3.2}
\rho(G_1)\leq\rho(G_2),
\end{align}
where the equality occurs if and only if $G_1=G_2$. For $s=1$, we admit $G_2=K_1\vee(K_{n-4}\cup K_2\cup K_1)$. Using \eqref{eq:3.1} and \eqref{eq:3.2}, we infer
$$
\rho(G)\leq\rho(K_1\vee(K_{n-4}\cup K_2\cup K_1)),
$$
with equality if and only if $G=K_1\vee(K_{n-4}\cup K_2\cup K_1)$. Observe that $K_1\vee(K_{n-4}\cup K_2\cup K_1)$ contains no $H$-factor, a contradiction. We are to consider $s\geq2$

\noindent{\bf Case 1.} $n\geq2s+5$.

Recall that $G_2=K_s\vee(K_{n-2s-2}\cup K_2\cup sK_1)$. For the partition $V(G_2)=V(K_s)\cup V(K_{n-2s-2})\cup V(K_2)\cup V(sK_1)$, the quotient matrix of $A(G_2)$ can be written as
\begin{align*}
B_2=\left(
  \begin{array}{cccc}
  s-1 & n-2s-2 & 2 & s\\
  s & n-2s-3 & 0 & 0\\
  s & 0 & 1 & 0\\
  s & 0 & 0 & 0\\
  \end{array}
\right).
\end{align*}
Hence, the characteristic polynomial of $B_2$ equals
\begin{align*}
\varphi_{B_2}(x)=&x^{4}-(n-s-3)x^{3}-(s^{2}+2s+1)x^{2}+(s^{2}n+2sn+n-2s^{3}-6s^{2}-7s-3)x\\
&-s^{2}n+2s^{3}+3s^{2}.
\end{align*}
One easily see that the partition $V(G_2)=V(K_s)\cup V(K_{n-2s-2})\cup V(K_2)\cup V(sK_1)$ is equitable. From Lemma 2.3, $\rho(G_2)$ is the largest root of $\varphi_{B_2}(x)=0$.

Let $G_*=K_1\vee(K_{n-4}\cup K_2\cup K_1)$. By the partition $V(G_*)=V(K_1)\cup V(K_{n-4})\cup V(K_2)\cup V(K_1)$, the quotient matrix of $A(G_*)$ is given by
\begin{align*}
B_*=\left(
  \begin{array}{cccc}
  0 & n-4 & 2 & 1\\
  1 & n-5 & 0 & 0\\
  1 & 0 & 1 & 0\\
  1 & 0 & 0 & 0\\
  \end{array}
\right),
\end{align*}
and so its characteristic polynomial is
$$
\varphi_{B_*}(x)=x^{4}-(n-4)x^{3}-4x^{2}+(4n-18)x-n+5.
$$
Since the partition $V(G_*)=V(K_1)\cup V(K_{n-4})\cup V(K_2)\cup V(K_1)$ is equitable, it follows from Lemma 2.3 that $\rho(G_*)$ is the largest root of $\varphi_{B_*}(x)=0$.

Note that $\varphi_{B_2}(\rho(G_2))=0$. Thus, we possess
\begin{align}\label{eq:3.3}
\varphi_{B_*}(\rho(G_2))=\varphi_{B_*}(\rho(G_2))-\varphi_{B_2}(\rho(G_2))=(s-1)f(\rho(G_2)),
\end{align}
where $f(\rho(G_2))=-\rho^{3}(G_2)+(s+3)\rho^{2}(G_2)-((s+3)n-2s^{2}-8s-15)\rho(G_2)+(s+1)n-2s^{2}-5s-5$.

Since $K_{n-s-2}$ is a proper subgraph of $G_2=K_s\vee(K_{n-2s-2}\cup K_2\cup sK_1)$, it follows from Lemma 2.1 that
\begin{align}\label{eq:3.4}
\rho(G_2)>\rho(K_{n-s-2})=n-s-3.
\end{align}

Let $f(x)=-x^{3}+(s+3)x^{2}-((s+3)n-2s^{2}-8s-15)x+(s+1)n-2s^{2}-5s-5$, $x\geq n-s-3$. Then the derivative function of $f(x)$ is
$$
f'(x)=-3x^{2}+2(s+3)x-(s+3)n+2s^{2}+8s+15, \ \ \ \ \ x\geq n-s-3.
$$
Notice that $\frac{s+3}{3}<n-s-3\leq x$ by $s\geq2$ and $n\geq2s+5$. Thus, we infer
\begin{align*}
f'(x)\leq&f'(n-s-3)\\
=&-3n^{2}+7(s+3)n-3s^{2}-22s-30\\
\leq&-3(2s+5)^{2}+7(s+3)(2s+5)-3s^{2}-22s-30 \ \ \ \ \ \ \ (\mbox{since} \ n\geq2s+5)\\
=&-s^{2}-5s\\
<&0.
\end{align*}
This implies that $f(x)$ is decreasing with respect to $x\geq n-s-3$. Combining this with \eqref{eq:3.4}, we deduce
\begin{align}\label{eq:3.5}
f(\rho(G_2))<&f(n-s-3)\nonumber\\
=&-(n-s-3)^{3}+(s+3)(n-s-3)^{2}\nonumber\\
&-((s+3)n-2s^{2}-8s-15)(n-s-3)+(s+1)n-2s^{2}-5s-5\nonumber\\
=&-n^{3}+(3s+9)n^{2}-(2s^{2}+15s+20)n+2s^{2}+10s+4.
\end{align}
Let $g(n)=-n^{3}+(3s+9)n^{2}-(2s^{2}+15s+20)n+2s^{2}+10s+4$. Using a similar discussion as above, we easily see that $g(n)$ is decreasing with respect to $n\geq2s+5$. In terms of $s\geq2$
and $n\geq2s+5$, we get
\begin{align*}
g(n)\leq&g(2s+5)\\
=&-(2s+5)^{3}+(3s+9)(2s+5)^{2}-(2s^{2}+15s+20)(2s+5)+2s^{2}+10s+4\\
=&-2s^{2}+4\\
<&0.
\end{align*}
Combining this with \eqref{eq:3.3}, \eqref{eq:3.5} and $s\geq2$, we obtain
$$
\varphi_{B_*}(\rho(G_2))=(s-1)f(\rho(G_2))<(s-1)g(n)<0,
$$
which yields that
\begin{align}\label{eq:3.6}
\rho(G_2)<\rho(G_*).
\end{align}
By virtue of \eqref{eq:3.1}, \eqref{eq:3.2} and \eqref{eq:3.6}, we conclude
$$
\rho(G)\leq\rho(G_1)\leq\rho(G_2)<\rho(G_*)=\rho(K_1\vee(K_{n-4}\cup K_2\cup K_1)),
$$
which contradicts $\rho(G)\geq\rho(K_1\vee(K_{n-4}\cup K_2\cup K_1))$.

\noindent{\bf Case 2.} $n=2s+4$.

In this case, $G_2=K_s\vee(2K_2\cup sK_1)$. Consider the partition $V(G_2)=V(K_s)\cup V(2K_2)\cup V(sK_1)$, the corresponding quotient matrix of $A(G_2)$ is
\begin{align*}
B_3=\left(
  \begin{array}{ccc}
    s-1 & 4 & s\\
    s & 1 & 0\\
    s & 0 & 0\\
  \end{array}
\right).
\end{align*}
The characteristic polynomial of $B_3$ equals
$$
\varphi_{B_3}(x)=x^{3}-sx^{2}-(s^{2}+3s+1)x+s^{2}.
$$
Note that the partition $V(G_2)=V(K_s)\cup V(2K_2)\cup V(sK_1)$ is equitable. By Lemma 2.3, $\rho(G_2)$ is the largest root of $\varphi_{B_3}(x)=0$. That is, $\varphi_{B_3}(\rho(G_2))=0$.

If $s\in\{2,3\}$, then $n\in\{8,10\}$, which contradicts $n\geq11$. Next, we discuss $s\geq4$.

Since $K_s\vee(s+4)K_1$ is a proper subgraph of $G_2=K_s\vee(2K_2\cup sK_1)$, it follows from Lemma 2.1 that
\begin{align}\label{eq:3.7}
\rho(G_2)>\rho(K_s\vee(s+4)K_1).
\end{align}
In terms of the partition $V(K_s\vee(s+4)K_1)=V(K_s)\cup V((s+4)K_1)$, the quotient matrix of $A(K_s\vee(s+4)K_1)$ equals
\begin{align*}
B_4=\left(
  \begin{array}{ccc}
    s-1 & s+4\\
    s & 0 \\
  \end{array}
\right).
\end{align*}
The characteristic polynomial of $B_4$ is
$$
\varphi_{B_4}(x)=x^{2}-(s-1)x-s(s+4).
$$
Note that the partition $V(K_s\vee(s+4)K_1)=V(K_s)\cup V((s+4)K_1)$ is equitable. According to Lemma 2.3, $\rho(K_s\vee(s+4)K_1)$ is the largest root of $\varphi_{B_4}(x)=0$. Together with \eqref{eq:3.7} and
$s\geq2$, we conclude
\begin{align}\label{eq:3.8}
\rho(G_2)>\rho(K_s\vee(s+4)K_1)=\frac{s-1+\sqrt{5s^{2}+14s+1}}{2}\geq\frac{3}{2}s+1.
\end{align}

Recall that $n=2s+4$, $\varphi_{B_*}(x)=x^{4}-(n-4)x^{3}-4x^{2}+(4n-18)x-n+5$, $\varphi_{B_3}(x)=x^{3}-sx^{2}-(s^{2}+3s+1)x+s^{2}$ and $\varphi_{B_3}(\rho(G_2))=0$. Thus, we have
\begin{align}\label{eq:3.9}
\varphi_{B_*}(\rho(G_2))=&\varphi_{B_*}(\rho(G_2))-\rho(G_2)\varphi_{B_3}(\rho(G_2))\nonumber\\
=&-s\rho^{3}(G_2)+(s^{2}+3s-3)\rho^{2}(G_2)-(s^{2}-8s+2)\rho(G_2)-2s+1.
\end{align}
Let $h(x)=-sx^{3}+(s^{2}+3s-3)x^{2}-(s^{2}-8s+2)x-2s+1$, $x\geq\frac{3}{2}s+1$. Then the derivative function of $h(x)$ is
$$
h'(x)=-3sx^{2}+2(s^{2}+3s-3)x-s^{2}+8s-2, \ \ \ \ \ x\geq\frac{3}{2}s+1.
$$
Since $s\geq4$ and $\frac{2(s^{2}+3s-3)}{6s}<\frac{3}{2}s+1\leq x$, we possess
\begin{align*}
h'(x)\leq&h'\Big(\frac{3}{2}s+1\Big)\\
=&-3s\Big(\frac{3}{2}s+1\Big)^{2}+2(s^{2}+3s-3)\Big(\frac{3}{2}s+1\Big)-s^{2}+8s-2\\
=&-\frac{15}{4}s^{3}+s^{2}-\frac{5}{2}s-11\\
<&0,
\end{align*}
which implies that $h(x)$ is decreasing with respect to $x\geq\frac{3}{2}s+1$. By \eqref{eq:3.8} and $s\geq4$, we obtain
\begin{align*}
h(\rho(G_2))<&h\Big(\frac{s-1+\sqrt{5s^{2}+14s+1}}{2}\Big)\\
=&-s\Big(\frac{s-1+\sqrt{5s^{2}+14s+1}}{2}\Big)^{3}+(s^{2}+3s-3)\Big(\frac{s-1+\sqrt{5s^{2}+14s+1}}{2}\Big)^{2}\\
&-(s^{2}-8s+2)\Big(\frac{s-1+\sqrt{5s^{2}+14s+1}}{2}\Big)-2s+1\\
=&\frac{1}{2}(-s^{4}+8s^{3}+28s^{2}-28s+3-(s^{3}+s^{2}-s-1)\sqrt{5s^{2}+14s+1})-\frac{1}{4}\\
<&\frac{1}{2}(-s^{4}+8s^{3}+28s^{2}-28s+3-(s^{3}+s^{2}-s-1)(2s+3))\\
=&\frac{1}{2}(-3s^{4}+3s^{3}+27s^{2}-23s+6)\\
<&0.
\end{align*}
Combining this with \eqref{eq:3.9}, we deduce
$$
\varphi_{B_*}(\rho(G_2))=h(\rho(G_2))<0,
$$
which leads to $\rho(G_2)<\rho(G_*)$. Together with \eqref{eq:3.1} and \eqref{eq:3.2}, we conclude
$$
\rho(G)\leq\rho(G_1)\leq\rho(G_2)<\rho(G_*)=\rho(K_1\vee(K_{n-4}\cup K_2\cup K_1)),
$$
which contradicts $\rho(G)\geq\rho(K_1\vee(K_{n-4}\cup K_2\cup K_1))$. Theorem 1.2 is proved. \hfill $\Box$

\medskip

\section*{Data availability statement}

My manuscript has no associated data.

\section*{Declaration of competing interest}

The authors declare that they have no conflicts of interest to this work.

\section*{Acknowledgments}

This work was supported by the Natural Science Foundation of Jiangsu Province (Grant No. BK20241949). Project ZR2023MA078 supported by Shandong
Provincial Natural Science Foundation.

\end{document}